Cubic equations of Babylonian mathematics

Kazuo MUROI

§1. Introduction

In Babylonian mathematics, cubic equations occur in several problems mainly concerning digging a hole, which must be the oldest known examples of the cubic equation. Naturally it was impossible for the Babylonian scribes of the Old Babylonian period (2000-1600 BCE) to reach the solution of the general cubic equation, because they did not develop appropriate mathematical symbols, nor did they conceive negative numbers and irrational numbers, not to mention complex numbers. Nevertheless, using a table of cube roots or the like they attempted to solve certain types of the cubic equations, of course including the simplest one, that is, $x^3 = a$. If we analyze in detail these cubic equations, some of which have not been completely understood yet, we will realize the ingenuity or tour de force of the Babylonian scribes in solving the equations.

In this paper, I shall clarify three cubic equations of Babylonian mathematics in particular, whose solutions have not been fully explained; BM 85200, no.6 and no.7, and YBC 4669, B2.[1]

§2. Simple cubic equations

In order to show a typical example of the simple Babylonian cubic equations, I cite the tablet IM 54478,[2] on which only one problem is syllabically written in Akkadian. Before reading it, we had better confirm the metrological units employed in the text. In Babylonia as well as in Sumer, the length or the width of a land was usually measured in nindan (1 nindan ≈ 6m), while the height or depth in kùš (1 kùš = 1/12 nindan ≈ 50cm). Therefore, the area unit was (area)-sar (1 sar = 1 nindan$^2$ ≈ 36m$^2$), and the volume unit was also (volume)-sar (1 sar = 1 nindan$^2$·kùš ≈ 18m$^3$). If the depth of a hole (z) is x nindan, z = 12x in kùš. As we will see below in §4, the number 12 was sometimes modified by the Sumerian word bal "change", and so I translate it into "conversion constant".

IM 54478 , transliteration and translation

Obverse

1. *šum-ma ki-a-am i-ša-al-ka um-ma šu-ú-ma*

2. *ma-la uš-ta-am-hi-ru ú-ša-pí-il-ma*

3. *mu-ša-ar ù zu-uz$_4$ mu-ša-ri*

4. *e-pí-ri a-su-uh ki-ia uš-tam-hi-ir*

5. *ki ma-ṣí ú-ša-pí-il*

6. *at-ta i-na e-pí-ši-ka*

7. [ 1,30 *ù* ] 12 *lu-pu-ut-ma i-gi* 12 *pu-ṭú-ur-ma*

8. [ 5 *ta-mar a-na* 1,]30 *e-pí-ri-ka*

Reverse

1. *i-ši-ma* 7,30 *ta-mar* 7,30

2. *mi-nam* íb-si$_8$ 30 íb-si$_8$ 30 *a-na* 1

3. *i-ši-ma* 30 *ta-mar* 30 *a-na* 1 *ša-ni-im*

4. *i-ši-ma* 30 *ta-mar* 30 *a-na* 12

5. *i-ši-ma* 6 *ta-mar* 30 *mi-it-ha-ar-ta-ka*

6. 6 *šu-pu-ul-ka*

Obverse

1. If he asks you so, (he says) thus:

2. I excavated deep as much as (the side of) the square that I made, and

3,4. I removed one and a half *mūšarum* (= sar) of earth.

4,5. How did I make the square? How deep did I excavate?

6. When you perform the calculations,

7. inscribe [1;30 and] 12, and make the reciprocal of 12, and

8. [ you see 0;5. By 1;]30 of your earth,

Reverse

1. multiply (0;5), and you see 0;7,30.

2. What is the cube root of 0;7,30?   0;30 is the cube root.

2,3. Multiply 0;30 by 1, and you see 0;30.

3,4. Multiply 0;30 by another 1, and you see 0;30.

4,5. Multiply 0;30 by 12, and you see 6.

5,6. 0;30 is the side of your square. 6 is your depth.

### Mathematical commentary

In this problem, if we denote the side of a cubic hole by x, the following

cubic equation is presented, and solved step by step :

$$12x^3 = 1;30$$

$$x^3 = \overline{12} \cdot 1;30 = 0;5 \cdot 1;30 = 0;7,30$$

$$\therefore x = \sqrt[3]{0;7,30} = 0;30 \,(\text{nindan}) = 6\,(\text{kùš}).$$

It seems that the scribe of this tablet used a table of cube roots;

n³-e n ba-si₍₈₎ (n = 1,2,3, ···, 60)⁽³⁾,

which literally means "n³ corresponds to n", giving $\sqrt[3]{7,30,0} = 30$ for n = 30.

The same cubic equation occurs in problem no. 22 of BM 85200, though written in compact style.

§3. More difficult cubic equation

As far as I know, no one has classified the problem B2, according to Neugebauer's numbering, of the tablet YBC 4669 as a cubic equation, probably because they have not seized its mathematical meaning owing to the fact that the answer is not given in the text and there is a certain mistake in writing a number by the scribe. However, this short problem obviously deals with a cubic equation which is different from those deciphered so far.

YBC 4669, Reverse I, lines 8-10 (= problem B2)

8. uš *ù* sahar

9. gar-gar-*ma* 33,20 (sic)

10. uš sahar en-nam

8. The length and the volume,

9. I added together, and (the result is) 33,22(!).

10. What are the length (and) the volume?

Mathematical commentary

Since the equation $12x^3 + x = 33,20$ does not have a rational solution, the cubic equation originally intended by the scribe would be;

$12x^3 + x = 33,22$.

Multiplying both sides by 18, we obtain

$(6x)^3 + 3(6x) = 33,22 \cdot 18 = 10,0,36$.

Successively calculating the values of $n^3 + 3n$ (n = 1, 2, 3, ···), we may notice

$33^3 + 3 \cdot 33 = 9,58,57 + 1,39 = 10,0,36$.

Therefore, $6x = 33$ and $x = 5;30$, though it is not given in the text.

The scribe who had made this problem may have known the fact that $x = m + (1/2)$ (m: integer) gives whole numbers to the cubic polynomial $12x^3 + x$ ;

$$12x^3 + x = 12\left(m + \frac{1}{2}\right)^3 + m + \frac{1}{2}$$

$$= 12m^3 + 18m^2 + 10m + \mathbf{2}.$$

If so, there is a strong possibility that he had forgotten to add 2 to $12 \cdot 5^3 + 18 \cdot 5^2 + 10 \cdot 5 = 2000 = 33,20$ at the last step of making the equation.

At any rate it is apparent that mathematicians of, say, the 16th century BCE were unable to solve the general cubic equation $x^3 + ax = b$, because

Italian mathematicians of the 16th century finally solved it after great struggle.

§4. Simultaneous cubic equations

There occur three simultaneous cubic equations in the problems 5, 6, and 7 of partly broken tablet BM 85200, of which no.5 is particularly well-known as well as no.23 both of which deal with the same cubic equation. The reconstructed mathematical meaning of no.5 is, according to Neugebauer, as follows:

$$xy + xyz = 1;10 \quad y = 0;40x \quad z = 12x,$$

from which by elimination of y and z, we obtain

$$8x^3 + 0;40x = 1;10.$$

Multiplying both sides by $12^2/0;40$, the equation becomes

$$(12x)^3 + (12x)^2 = z^3 + z^2 = 4,12.$$

Since we have a table of such as:

$n^3 + n^2$ -e n ba-si (n =1, 2, 3, ⋯,60)[4], which literally means "$n^3 + n^2$ corresponds to n",

we can easily obtain the solution $12x = 6$, that is, $x = 0;30$. Although there seems to be no question about this interpretation which is believed to be right by most of historians of mathematics, we cannot explain the next problems nos.6 and 7 along the same line. In my judgment, the key to solving the above three problems of BM 85200 is not the use of some tables but the application of a Babylonian mathematical method called *makṣarum* "factorization". In fact, this term occurs in one problem dealing with a cubic equation[5]:

> *[ma]-ak-ṣa-ru-um ša* ba-si (YBC 6295, line 1) "the factorization method of a cube root".

Now let us analyze the problem no.7 first, because we can clearly understand its "factorization method" for the sake of its given data that leads to simple simultaneous cubic equations.

   BM 85200, Obverse I, lines 15-20( = no.7 ),transliteration and translation

15.  túl-sag *ma-la* uš GAM-*ma* 1 sahar-hi-a ba-zi [k]i-*ri ù* sahar-hi-a ul-<gar> 1,10 uš ugu sag 10 dirig

16.  za-e 1 *ù* 12 [b]al gar-ra 10 [dirig *a-n*]*a* 1 *i-ši* 10 *ta-mar a-na* 12 *i-ši* 2 *ta-mar*

17.  10 *šu-tam-<hir>* 1,40 *ta-mar a-na* 2 *i-ši* 3,2[0 *t*]*a-mar* igi 3,20 du₈-*a* 18 *ta-mar*

18.  *a-na* 1,10 *i-ši* 21 *ta-mar*  3  2  21(sic)  íb-si₈ [ 10 *a-na* 3 *i*]*-ši* 30 uš

19.  10 *a-na* 2 *i-ši* 20 sag 3 *a-na* 2 *i-ši* 6 *ta-mar* [6] GAM

20.  *ne-pé-šu*[*m*]

15.  A very old well. The depth is equal to the length. I removed the earth, 1(sar). I added the area and the volume, (and) 1;10(sar). The length exceeds the width by 0;10(nindan).

16. You, put down 1 and 12, the conversion constants. Multiply 0;10, the excess, by 1, (and) you see 0;10. Multiply (0;10,the excess,) by 12, (and) you see 2.

17. Square 0;10, (and) you see 0;1,40. Multiply (the result) by 2, (and) you see 0;3,20. Make the reciprocal of 0;3,20, (and) you see 18.

18. Multiply (the result) by 1;10, (and) you see 21.   3   2   3;30(!) are the roots. Multiply 0;10 by 3, (and) 0;30 is the length.

19. Multiply 0;10 by 2, (and) 0;20 is the width. Multiply 3 by 2, (and) you see 6. 6 is the depth.

20. The procedure.

## Mathematical commentary

x : the length,   y : the width,   z : the depth.

z = 12x ---(1)   xyz = 1 --- (2) xy + xyz = 1;10 --- (3) x − y = 0;10 --- (4)

From (1) and (2), we obtain $12x^2y = 1$ --- (5)

Dividing both sides of (4) by 0;10, the text gives;

$$\frac{x}{0;10} - \frac{y}{0;10} = X - Y = 1 \quad \text{---} \quad (6)$$

Dividing both sides of (5) by 12(0;10)³ , the text also gives;

$$\left(\frac{x}{0;10}\right)^2 \left(\frac{y}{0;10}\right) = \frac{1}{0;10^2} \cdot \frac{1}{0;10 \cdot 12} = \frac{1}{0;1,40} \cdot \frac{1}{2}$$

$$= \frac{1}{0;3,20} = 18, \text{ that is,}$$

$X^2Y = 18 \ ( = 3^2 \cdot 2) \ \cdots \ (7)$

On the other hand, from (3) we obtain

$xy(z + 1) = 1;10$, which is transformed into (8) by the same division;

$$\frac{x}{0;10} \cdot \frac{y}{0;10} \cdot \frac{z+1}{0;10 \cdot 12} = 1;10 \cdot 18 = 21, \text{ or}$$

$$XY \cdot \frac{z+1}{2} = 21 \ \cdots \ (8)$$

Observing (6) and (7), we easily notice the solution

$X = 3$ and $Y = 2$, resulting in $(z + 1)/2 = 3;30$ from (8).

These three numbers 3, 2, and 3;30 are called íb-si$_8$ "roots" in line 18. Therefore,

$x = 0;10 \cdot 3 = 0;30$, $y = 0;10 \cdot 2 = 0;20$, and from $z/2 + 0;30 = 3;30$, we finally

obtain the depth, $z = 3 \cdot 2 = 6$.

Thus the Babylonian scribe factorized 18 into $3^2 \cdot 2$ and found the sole solution (X, Y) = (3, 2), from which x, y, and z were immediately obtained.

BM 85200, Obverse I, lines 9-14 ( = no.6 ), transliteration and translation

9. túl-sag *ma-la* uš GAM-*ma* 1 sahar-hi-a ba-zi ki-*ri ù* sahar-hi-a ul-gar 1,10 uš *ù* sag 50 uš sag en-\<nam\>

10. za-e 50 *a-na* 1 bal *i-ši* 50 *ta-mar* 50 *a-na* 12 *i-ši* 10 *ta-mar*

11. 50 *šu-tam-*\<hir\> 41,40 *ta-mar a-na* 10 *i-ši* 6,56,40 *ta-mar* igi-*šu* du$_8$-a 8,38,24 *ta-*\<mar\>

12.  *a-na* 1,10 *i-ši* 10,4,48 *ta-mar*   36   24   42   íb-si₈

13.  36 *a-na* 50 *i-ši* 30 uš 24 *a-na* 50 *i-ši* 20 sag 36 *a-na* 10 <*i-ši*> 6 GAM

14.  [*n*]*e-pé-šum*

9.  A very old well. The depth is equal to the length. I removed the earth, 1(sar). I added the area and the volume, (and) 1;10(sar).( The sum of ) the length and the width is 0;50(nindan). What are the length (and) the width?

10. You, multiply 0;50 by the conversion constant 1, (and) you see 0;50. Multiply 0;50 by (the conversion constant) 12, (and) you see 10.

11. Square 0;50, (and) you see 0;41,40. Multiply (the result) by 10, (and) you see 6;56,40. Make its reciprocal,(and) you see 0;8,38,24.

12. Multiply (the result) by 1;10, (and) you see 0;10,4,48.   0;36   0;24   0;42   are the roots.

13. Multiply 0;36 by 0;50, (and) 0;30 is the length. Multiply 0;24 by 0;50, (and) 0;20 is the width. Multiply 0;36 by 10, (and) 6 is the depth.

14. The procedure.

<div style="text-align:center">Mathematical commentary</div>

$z = 12x$ ---(1)   $xyz = 1$ ---(2)   $xy + xyz = 1;10$ ---(3)   $x + y = 0;50$ ---(4)

From (1) and (2), we obtain $12x^2 y = 1$ ---(5)

Dividing both sides of (4) by 0;50, the text gives

$$\frac{x}{0;50} + \frac{y}{0;50} = X + Y = 1 \quad \text{---(6)}$$

Similarly dividing both sides of (5) by $12(0;50)^3$, we obtain;

$$\left(\frac{x}{0;50}\right)^2\left(\frac{y}{0;50}\right) = \frac{1}{0;50^2} \cdot \frac{1}{0;50 \cdot 12}$$

$$= \frac{1}{0;41,40} \cdot \frac{1}{10} = \frac{1}{6;56,40}$$

$$= 0;8,38,24, \text{ that is,}$$

$X^2Y = 8,38,24 \cdot 60^{-3} = 2^7 \cdot 3^5 \cdot 60^{-3}$

$= 72^2 \cdot 6 \cdot 60^{-3}$ or $36^2 \cdot 24 \cdot 60^{-3}$

$= (1;12)^2 (0;6)$ or $(0;36)^2(0;24)$ --- (7) ( because $X \geq Y$)

Observing (6) and (7) we may notice $X = 0;36$ and $Y = 0;24$. On the other hand, from (3) we also obtain;

$xy(z + 1) = 1;10$, which is transformed into (8) using the results just calculated above;

$$\frac{x}{0;50} \cdot \frac{y}{0;50} \cdot \frac{z+1}{0;50 \cdot 12} = 1;10 \cdot 0;8,38,24$$

$= 10,4,48 \cdot 60^{-3} = 2^6 \cdot 3^4 \cdot 7 \cdot 60^{-3}$

$= 36 \cdot 24 \cdot 42 \cdot 60^{-3}$

$= 0;36 \cdot 0;24 \cdot 0;42$, or

$$XY \cdot \frac{z+1}{10} = 0;36 \cdot 0;24 \cdot 0;42 \quad \text{---} \quad (8).$$

Therefore,

x = 0;36 · 0;50 = 0;30, y = 0;24 · 0;50 = 0;20, and from (8), z/10 + 0;6 = 0;42,

we finally obtain the depth z = 0;36 · 10 = 6.

Again the Babylonian scribe factorized the number 8,38,24 into $36^2 \cdot 24$, and found one of the solutions ( X, Y ) = ( 0;36, 0;24 ), from which x, y, and z were obtained in the same way. As to another solution of (6) and (7);

(X, Y) = ( $(12 + 12\sqrt{7})/60$, $(48 - 12\sqrt{7})/60$ ),

it would have been impossible for the Babylonians to get it.

Considering the above two problems, it is more likely that the scribe of our tablet have used the same method, the factorization of a number, in solving the problem no.5:

$z^3 + z^2 = z^2( z + 1 ) = 4,12 = 6^2 \cdot 7 = 6^2(6 + 1)$ ∴ z = 6.

In other words, he only found the solution not referring to the table "$n^3 + n^2$-e n ba-si", mentioned above.

The Babylonian factorization method may be useful for searching one positive rational solution, if any, of a numerical cubic equation of the type $x^3 + ax^2 = b$, but it would be almost useless for $x^3 + ax = b$, of which we have only one example discussed in §3 above. For the same method probably used in Chinese mathematics of the 7th century, see Appendix below.

Appendix

In his work of Jigu Suanjing, a Chinese mathematician Wang Xiaotong,

flourished in the first half of the 7th century, presented a system of simultaneous equations;

$$x^2 + y^2 = z^2, \quad xy = 706\frac{1}{50}(=P), \quad z - x = 36\frac{9}{10}(=S).$$

According to his directions, we get the following cubic equation after the elimination of y and z;

$$x^3 + \frac{S}{2}x^2 = \frac{P^2}{2S} \quad .(7)$$

Since he gave the answer x = 14 + (7/20) without any specific explanation, we must conjecture his method of solution. If we factorize S/2 and P²/2S, we get;

$$x^3 + \frac{3^2 \cdot 41}{2^2 \cdot 5}x^2 = \frac{7^2 \cdot 41^3}{2^2 \cdot 5^3}$$

This clearly suggests to us;

$$x = \frac{41n}{2^2 \cdot 5} \quad (n: \text{positive integer}).$$

Therefore, the above cubic equation reduces to;

n³ + 9n² = n²(n + 9) = 7²(7 + 9),

from which n = 7 is immediately obtained. Thus Wang had got the answer x = 14 + (7/20), y = 49 + (1/5), and z = 50 + (1/4), I believe.

Notes

(1) O. Neugebauer, *Mathematische Keilschrift-Texte* (= MKT) Ⅰ, Ⅱ (1935), Ⅲ(1937). BM 85200, no.6 and no.7; MKT I, pp.193-219. YBC 4669, B2; MKT Ⅲ, p.64.

(2) T. Baqir, Some more mathematical texts from Tell Harmal, *Sumer* 7 (1951), pp. 28-45.

(3) MKT Ⅰ, p. 73.

(4) MKT Ⅰ, pp. 76-77.

(5) K. Muroi, Extraction of Cube Roots in Babylonian Mathematics, *Centaurus*, vol. 31(1989), pp. 181-188. For the meaning of *makṣarum*, see: E. M. Bruins, Some Mathematical Texts, *Sumer* 10 (1954), pp. 55-61, esp. 56, and J. Friberg, *A Survey of Publications on Sumero-Akkadian Mathematics, Metrology and Related Matters (1854-1982)*, p.64.

(6) As Neugebauer pointed out, the system of simultaneous cubic equations of no.7 can be reduced to the following simultaneous quadratic equations:

$$z = 12x, \quad xy + 1 = 1;10, \quad x - y = 0;10.$$

But the scribe of our tablet did not do so.

(7) Y. Mikami, *Development of Mathematics in China and Japan* (1913), p.54.

Acknowledgement

I am grateful to Prof. S. Nakamura, who had read my manuscript and gave me many helpful suggestions.